\documentclass[10pt]{article}

\newcommand{\qed}{\hfill \quad \rule{1.0 mm}{3.0 mm}}

\newcommand{\subterm}{\sqsubseteq}
\newcommand{\rotmap}{{\rm rot}}
\newcommand{\refmap}{{\rm ref}}

\begin{document}

\begin{center}
{\large \bf A squarefree term not occurring in the Leech sequence}
\end{center}

\vspace{2ex}
\begin{center}
Benjamin Wells \\ 
University of San Francisco \\
March 22, 2020
\end{center}

\bigskip

\noindent
{\bf Introduction}

\vspace{0.1in}

Let $T = \{A,B,C\}$ be a finite alphabet, and let $P$ be the substitution map from $T^+$ to $T^+$ (a semigroup homomorphism) determined by
\[
\begin{array}{c}
P(A) = \overline{A} =  ABCBA\ CBC\ ABCBA, \\
P(B) = \overline{B} =  BCACB\ ACA\ BCACB, \\
P(C) = \overline{C} =  CABAC\ BAB\ CABAC. \\
\end{array}
\]
The Leech sequence $L$ \cite{L} is the squarefree sequence obtained as the limit of the palindromes
\[
A, \overline{A}, \overline{\overline{A}}, \ldots, P^n(A), \ldots.
\]
In order to specify a certain class of pseudorecursive varieties of semigroups \cite{W}, it is helpful to have a squarefree term in 3 variables such that no substitution instance occurs as a subterm of $L$.  We show that
\[
\kappa_1 = aba\ cbc\ aba\ c
\]
is such a term.  Except for one case, the term
\[
\kappa_2 = aba\ cbc\ aba
\]
will serve, and we focus on it.  If $\ \hat{} : \{a,b,c\} \rightarrow T^+$, then we write $\langle\hat{a},\hat{b},\hat{c}\rangle = \hat{a}\hat{b}\hat{a}\ \hat{c}\hat{b}\hat{c}\ \hat{a}\hat{b}\hat{a}$.  The terms $\hat{a},\hat{b},\hat{c}$ are called {\it keywords} of the map $\hat{}$\ .  Note that $|\langle \hat{a},\hat{b},\hat{c}\rangle| = 4|\hat{a}|+3|\hat{b}|+2|\hat{c}|$.

We can replace $\kappa_1$ with $\kappa_3 = c\ aba\ cbc\ aba$, and if a palindromic term is desired, we can take  $\kappa_4 = c\ aba\ cbc\ aba\ c$.

\bigskip
 
\noindent
{\bf Rigid and flush keywords.}

An $L$-{\it block} is $\overline{A}$, $\overline{B}$, or $\overline{C}$. A {\it matrix} $\mu$ of a subterm $\alpha \subterm L$ is a minimal-length succession of $L$-blocks such that $\alpha \subterm \mu$. Thus, there is a (unique) term $\nu \in T^+$ such that $\overline{\nu} = \mu$.  It is not immediately obvious that $\nu \subterm L$, but it becomes so when we choose $n \geq1$ such that $\mu \subterm P^n(A)$, for then $\nu \subterm P^{n-1}(A)$.  In general, the term $\alpha$ can have several unequal matrices and can be located in more than one position in a given matrix.  For example, $BACAB \subterm \overline{B},\overline{AC},\overline{CA}$, and $ABCBA \subterm \overline{A}$ only, but in two positions:  $\underline{ABCBA}\ CBC\, \underline{ABCBA}$. 

We write $\alpha\ (\mu_1,\mu_2, \ldots, \mu_j)$ to indicate the matrices of the term $\alpha \subterm L$.  Thus, $BACAB \ (\overline{B},\overline{AC},\overline{CA})$ and $ABCBA\ (\overline{A},\overline{A})$.  A subterm of $L$ is {\it rigid} iff it has a unique matrix and appears in only one position in it.  A rigid term is {\it left flush} if it is an initial segment of its matrix; it is {\it right flush} if it is a final segment of its matrix.  If it is flush on both sides (and hence a string of blocks), then it is {\it blush}.  We can also use the terms flush and blush with respect to a particular matrix for a term, even if it has other matrices or other positions in a matrix; for example, $ABCBA$ is left and right flush in the matrix $\overline{A}$; because it isn't rigid, it isn't blush.

The {\it rotation map} rot on $T$ is given by:
\[
\begin{array}{c}
\rotmap(A) = B,\\
\rotmap(B) = C,\\
\rotmap(C) = A,
\end{array}
\]
naturally extended to a semigroup homomorphism on $T^+$. The unambiguous $\rotmap(L)$ is seen to be the limit of 
\[
B, \overline{B}, \overline{\overline{B}}, \ldots, P^n(B), \ldots.
\]
If $\alpha \subterm L$, then not quite obviously $\rotmap(\alpha) \subterm L$ by the observation that given $n>1$ such that $\rotmap(\alpha) \subterm P^n(B)$, then $\rotmap(\alpha) \subterm P^{n+1}(A)$.

The {\it reflection map} ref on $T$ is given by:
\[
\begin{array}{c}
\refmap(A) = A,\\
\refmap(B) = C,\\
\refmap(C) = B,
\end{array}
\]
which also extends to a semigroup homomorphism on $T^+$.  

Assume $\alpha\subterm \beta\subterm L$.  Then $\alpha$ rigid means $\beta$ and $\rotmap(\alpha)$ are rigid.  Furthermore $\alpha \subterm \overline{\alpha} \subterm \overline{\beta} \subterm L$, and $\rotmap(\alpha) \subterm \rotmap(\beta) \subterm L$. We have $\refmap(\alpha) \subterm \refmap(\beta) \subterm \refmap(L)$, but unlike the rotational symmetry observed, there is no reflection symmetry:  $\refmap(\alpha)$ is not necessarily a subterm of $L$.  In particular, $\refmap(\overline{A}), \refmap(\overline{B}), \refmap(\overline{C}) \not \subterm L$.  Indeed, $ABACABCB$ occurs in $\overline{CA}$, but $\refmap(ABACABCB) = ACABACBC \not \subterm L$.  This is the shortest length for a term with this property.

An anomaly occurs with $ACBABCA$:  it is rigid, appearing only in $\overline{C}$, but its reflection $ABCACBA$ occurs with matrices $\overline{AB}$ and $\overline{BA}$.

A {\it local instance} $[\hat{a},\hat{b},\hat{c}, \rho]$ of $\kappa_2$ occurs in $L$ when $\rho$ is a term such that $\rho \langle \hat{a},\hat{b},\hat{c} \rangle$ is an initial segment of $L$.  A keyword for a local instance $[\hat{a},\hat{b},\hat{c},\rho]$ is {\it locally rigid} if within this local instance, the keyword has a unique place in a unique matrix.  This is construed to mean four instances of $\hat{a}$, three instances of $\hat{b}$, and two instances of $\hat{c}$ according to the pattern of $\kappa_2$.  It ignores such things as $\hat{a} \subterm \hat{b}$.  Observe that all rigid keywords are locally rigid in all local instances of $\kappa_2$.  The length of a local instance is given by $|[\hat{a},\hat{b},\hat{c}, \rho]| = |\langle \hat{a},\hat{b},\hat{c} \rangle| = 4|\hat{a}|+3|\hat{b}|+2|\hat{c}|$.

\bigskip

{\bf \S 1 Rigidity of 7-letter keywords.}

\medskip

By inspection, all fifty-four (eighteen up to rotation) 7-letter subwords of $L$ are rigid except for $ABCACBA \ (\overline{AB}, \overline{BA})$ and its rot variants.  Let us deal with this special case.

\medskip

{\bf Theorem 1 (Seven-letter keywords).}  {\it Suppose $[\hat{a},\hat{b},\hat{c},\rho]$ is a local instance of $\kappa_2$ where some keyword is not locally rigid and has length 7.  The only assignment possible is
\[
\langle B,C,ABCACBA \rangle \subterm L,
\]
plus the rot variants.  None of the keywords is locally rigid, in fact all nine instances of the keywords have different matrices or different positions in a matrix.}

\medskip

{\it Proof.}  (I)  We investigate $\langle \hat{a}, \hat{b}, \hat{c} \rangle = \langle ABCACBA, -, - \rangle$ where two adjacent local occurrences of $\hat{a} = ABCACBA$ (i.e., there are no intervening local occurrences of $\hat{a}$) are not rigid (their matrices or positions in a single matrix differ).  Because $\hat{a}$ is not locally rigid, such a pair exists.

(A)  The first case is  $ABCACBA x ABCACBA$ with matrix $\overline{AB}  \cdots \overline{BA}$ (the shortest example is when there is a single $\overline{B}$ between the two $\overline{A}$s ($\overline{ABA}$) and $x = C$).  

(1) So if this is an instance of $\hat{a} \hat{b} \hat{a}$, then $x = \hat{b}$ has the form $C$ or the form  $CABCA\cdots ACBAC$.  From the matrix $\overline{AB}\cdots \overline{BA}$, $\hat{c}$ must have the form $BCBA \cdots $ or $\cdots ABCB$, depending on whether this $\hat{a} \hat{b} \hat{a}$ comes from the beginning or end of $\kappa_2$.  Alternatively, $\hat{c}$ can be simply $B$.   

(a)  If the former, then either $\hat{b}\hat{c}$ has a square on $CB$, or $\hat{c}\hat{b}$ has a square on $BC$.

(b) If the latter, $\hat{c} = B$, then there are the two forms for $B$ to consider.

(i)  If $\hat{b} = C$, then
\[
\hat{a}\hat{b}\hat{a}\hat{c}\hat{b}\hat{c}\hat{a} =ABCACBA\ C\ ABCACBA\ B\ C\ B\ AB \cdots
\]
fails to be a segment of $L$ at the last $B$ listed.  

(ii) If $\hat{b}$ takes the longer form, then $\hat{b}\hat{c}$ has a square on $ACB$, and $\hat{c}\hat{b}$ has a square on $BCA$.

(2) If this is an instance of $\hat{a}\hat{c}\hat{b}\hat{c}\hat{a}$, then
\[
x = \hat{c}\hat{b}\hat{c} = CABCA\cdots ACBAC.
\]    
Because $\hat{c}$ starts and ends with $C$, $\hat{b}$ has to start and end with $B$; otherwise, there will be a square on $A$ in $\hat{a}\hat{b}$ or $\hat{b}\hat{a}$.  This means  $\hat{c}$ takes the form $CABCA \cdots ACBAC$ ($\hat{c}$ cannot be both $CA$ and $AC$).   But that gives a square on $BCA$ in $\hat{b}\hat{c}$ and a square on $ACB$ in $\hat{c}\hat{b}$.

\medskip

(B)  The second case is $ABCACBA x ABCACBA$ with matrix $\overline{BA} \cdots \overline{AB} $ (the shortest example is when there is a single $\overline{A}$ between the two $\overline{B}$s ($\overline{BAB}$) and $x = BCBA\ CBC\ ABCB$).

(1) So if this is an instance of $\hat{a} \hat{b} \hat{a}$, then $x = \hat{b}$ has the form  $BCBA\cdots ABCB$.  This means for $\hat{c}$ to avoid a square on $A$ with $\hat{a}$ and on $B$ with $\hat{b}$, it must have the form $C\cdots C$ (possibly a single $C$).  Unfortunately, that gives a square on $BC$ in $\hat{b}\hat{c}$ and on $CB$ in $\hat{c}\hat{b}$. 

(2) If this is an instance of $\hat{a}\hat{c}\hat{b}\hat{c}\hat{a}$, then the shortest $\hat{c}\hat{b}\hat{c}$ can be is 
\[
x =BCBA\ CBC\ ABCB.
\]

(a)  If $\hat{c} = B$, then $\hat{b} = CBA\cdots ABC$ which yields an impossible square on $CBA$ in $\hat{a}\hat{b}$ and on $ABC$ in $\hat{b}\hat{a}$.  

(b)  If $\hat{c} = BCB\cdots BCB$ (even if just $BCB$), then $\hat{b} = A\cdots A$, which gives an impossible square on $A$ in both $\hat{a}\hat{b}$ and $\hat{b}\hat{a}$.

\medskip

(II)  We turn to $\langle \hat{a}, \hat{b}, \hat{c} \rangle = \langle -, ABCACBA, - \rangle$ where two adjacent local occurrences of $\hat{b} = ABCACBA$ are not rigid (since $\hat{b}$ is not locally rigid, such a pair always exists).  

(A)  The first case is  $ABCACBA x ABCACBA$ with matrix $\overline{AB}  \cdots \overline{BA}$ (the shortest example is $\overline{ABA}$ and $x = C$).  

(1)  Consider $x = \hat{a}\hat{c} = CABC \cdots CBAC$, for this cannot be the one letter $C$.  No matter how $x$ is divided, $\hat{a} = C \cdots$ and $\hat{c} = \cdots C$, so $\hat{c}\hat{a}$ has a square on $C$.

(2)  The same argument works for $x = \hat{c}\hat{a}$.

(B)  The second case is $ABCACBA x ABCACBA$ with matrix $\overline{BA} \cdots \overline{AB}$ (the shortest example is $\overline{BAB}$ and $x = BCBA\ CBC\ ABCB$).  

(1)  Consider $x = \hat{a}\hat{c} = BCB \cdots BCB$.  No matter how $x$ is divided, $\hat{a} = B \cdots$ and $\hat{c} = \cdots B$, so $\hat{c}\hat{a}$ has a square on $B$.

(2)  The same argument works for $x = \hat{c}\hat{a}$.

\medskip

(III)  Finally, we face $\langle \hat{a}, \hat{b}, \hat{c} \rangle = \langle -, -, ABCACBA \rangle$ where the two adjacent local occurrences of $\hat{c} = ABCACBA$ are not rigid. 

(A)  The first case is  $ABCACBA x ABCACBA$ with matrix $\overline{AB} \cdots \overline{BA}$ (the shortest example is when there is just $\overline{ABA}$ and $x = C$).  

(1)  If $\hat{b} = C$, then $\hat{a}$ has to start and end with $B$ to avoid squares, so there are two cases.

(a)  If $\hat{a} = B$, then in fact
\[
\langle \hat{a}, \hat{b}, \hat{c} \rangle = \langle B, C, ABCACBA \rangle = B\,C\,B\, ABCACBA\ C\, ABCACBA\ B\,C\, B \subterm L
\]
with matrix $\overline{ABA}$.  This is the single situation apart from its rot variants and $P$ expansions that admits instances of $\kappa_2$ in $L$.  We note that all three keywords are not locally rigid; indeed, all nine of the instances of the keywords have different matrices or positions.

(b)  If $\hat{a} = BCB \cdots BCB$ (even if just $BCB$), then $\hat{a}\hat{b}$ has a square on $BC$, and $\hat{b}\hat{a}$ has a square on $CB$.

(2)  If $x = \hat{b}$ takes the form $CABCA \cdots ACBAC$, then $\hat{a} = B \cdots B$ to avoid squares at the end letters. But that gives a square on $BCA$ in $\hat{a}\hat{b}$ and a square on $CBA$ in $\hat{b}\hat{a}$.

(B)  The second case is $ABCACBA x ABCACBA$ with matrix $\overline{BA} \cdots \overline{AB}$ (the shortest example is when there is just $\overline{BAB}$ and $x = BCBA\ \,CBC\ ABCB$).  The general (and only) form for $x$ is $\hat{b} = BCB \cdots BCB$.  That means $\hat{a} = C\cdots C$, possibly just $C$, to avoid squares at end letters.  But now $\hat{a}\hat{b}$ has a square on $CB$, and $\hat{b}\hat{a}$ has a square on $BC$.    \qed

\bigskip

The two 8-letter right extensions 
\[
ABCACBAB \mbox{\ and\ } ABCACBAC
\]
 of $ABCACBA$ appearing in $L$ are rigid, so this means all subwords of $L$ of seven letters or more are rigid, apart from this exceptional case.

\bigskip

{\bf \S 2 General rigidity results.}
\medskip

In this section, we examine how rigidity of a keyword constrains other keywords and restricts or prevents instances of $\kappa_2$ in $L$.

\bigskip

{\bf Theorem 2 (Locally rigid and flush keyword).}  {\it Suppose $[\hat{a},\hat{b},\hat{c},\rho]$ is a local instance of $\kappa_2$.  Assume no keyword is flush.  If one keyword is locally rigid, then the other two are rigid.  

In addition, there is an associated assignment $\,\check{}\,$ such that $|\check{a}| = |\hat{a}|, |\check{b}| = |\hat{b}|, |\check{c}| = |\hat{c}|$, and term $\rho'$ such that $[\check{a}, \check{b}, \check{c},\rho']$ is a local instance of $\kappa_2$, and one of the keywords $\check{a}, \check{b},\check{c}$ is locally rigid and flush.}

\medskip

{\it Proof.}  
(1)  We start by assuming $\hat{a}$ locally rigid in this local instance of $\kappa_2$.  Let its matrix be $M^a_1 M^a_2 \cdots M^a_j$.  Consider the left instance of $\hat{a}\hat{b}\hat{a}$, where the matrix of $\hat{b}$ is $M^{b,1}_1 M^{b,1}_2 \cdots M^{b,1}_k$.  Because $\hat{a}$ is not right flush, $\hat{b}$ will start in $M^a_j$, that is, $M^{b,1}_1 = M^a_j$.  Again, because $\hat{a}$ is not left flush, $\hat{b}$ will end in $M^a_1$, so $M^{b,1}_k = M^a_1$.  If $\hat{b}$ is a subterm of  $M^a_j M^a_1$ (or simply of $M^a_1$, if $j = 1$), $L$ contains an impossible square on $M^a_1 \cdots M^a_j$.  So $|\hat{b}| \geq 15$, meaning $\hat{b}$ is rigid.  If the leftmost $\hat{c}$ is a subterm of $M^a_j M^b_1$ (and these may coincide), then $L$ has a square.  Therefore $|\hat{c}| \geq 15$, and $\hat{c}$ is rigid.

(2)  Now assume $\hat{b}$ locally rigid in this local instance of $\kappa_2$; it is not flush.  If $\hat{a}\hat{c}$ is contained within $M^b_k M^b_1$ (or within $M^b_1$ if $k=1$), then there is an impossible square on $M^b_1 \cdots M^b_k$.  Therefore, $\hat{a}\hat{c}$ contains at least one block and has length $\geq 15$.  That means that either $|\hat{a}| \geq 8$ or $|\hat{c}| \geq 8$, so at least one of $\hat{a}$ or $\hat{c}$ is rigid.  

We start with $\hat{a}$ rigid and not flush.  So if the left $\hat{c}$ goes no farther right than $M^b_2$, there is an impossible square on $M^a_1 \cdots M^a_j M^b_2 \cdots M^b_k$ in $L$.  Therefore, this $\hat{c}$ spans a block, has length at least 13, and is rigid.

We turn to the alternative:  $\hat{c}$ rigid and not flush.  If $\hat{a}$ fits in $M^c_1$, then $L$ has the impossible square on $M^b_1\cdots M^b_k M^c_{1/2} \cdots M^c_l$.  Thus $|\hat{a}| \geq 13$, and $\hat{a}$ is rigid.  This completes the proof starting with $\hat{b}$.

(3)  Let's assume that $\hat{c}$ is locally rigid in this local instance of $\kappa_2$. The argument follows the proof for $\hat{a}$ above.

So in every case, two keywords are rigid and the third is at least locally rigid (recall none are flush), proving the first part of the theorem.

\medskip

We examine pieces of the underlying sequence $L$.  Let $\gamma$ be the initial segment of $M^a_1$ before $\hat{a}$, and let $\delta$ be the final segment of $M^a_j$ after $\hat{a}$.  Because $\hat{a}$ is not flush, $\gamma$ and $\delta$ are not empty.  We note that $\delta$ is an initial segment of $\hat{b}, \hat{c}$ (from $\hat{a}\hat{b}$ and $\hat{a}\hat{c}$), and $\gamma$ is a final segment of $\hat{b},\hat{c}$ (from $\hat{b}\hat{a}$ and $\hat{c}\hat{a}$).  

Let $m = \min(|\gamma|, |\delta|, 13-|\gamma|, 13-|\delta|)$.  Then $m$ is the shortest distance (in one direction or the other) from a keyword boundary to a block boundary.  Note that $m \leq 6$.  Shift all keyword boundaries in the correct direction by $m$ letters in $L$ to obtain the assignment $\,\check{}\,$ that has at least one locally rigid and flush keyword.  Coherence of keywords (all three instances of $\check{b}$, say, are equal) is guaranteed, because no block boundaries are crossed; in the same way, rigidity is preserved under this shift, proving more than sought.   Observe that $|\rho| \geq m$, since $\rho$ is at least as long as $\gamma$; $\rho'$ results from $\rho$ by adding or trimming $m$ letters---what they are is no matter.     \qed

\bigskip

{\bf Theorem 3 (Rigidity and blocks). }  {\it Suppose $[\hat{a},\hat{b},\hat{c},\rho]$ is a local instance of $\kappa_2$.  If one keyword is locally rigid and flush, then all three are rigid and blush.}

\medskip

{\it Proof.}  We start by assuming $\hat{a}$ locally rigid in this local instance of $\kappa_2$ and left flush.  Let its matrix be $M^a_1 M^a_2 \cdots M^a_j$, where $j = 1$ means there is only one block.  Consider the left instance of $\hat{a}\hat{b}\hat{a}$, where the matrix of $\hat{b}$ is $M^{b,1}_1 M^{b,1}_2 \cdots M^{b,1}_k$.  Because $\hat{a}$ is left flush, this $\hat{b}$ may start in $M^a_j$, that is, $M^{b,1}_1 = M^a_j$.  Because this $\hat{b}$ is right flush (ending at second instance of $\hat{a}$), either $\hat{b}$ is a final subterm of $M^a_j$ or ends with $M^{b,1}_k, k\geq 2$.  In the former case, $L$ has an impossible square on $M^a_1 \cdots M^a_j$ from $\hat{a}\hat{b}\hat{a}$.  In the latter, $|\hat{b}| \geq 13$ and is therefore rigid.

Now consider $\hat{a}\hat{c}\hat{b}$.  If $\hat{c} \subterm M^a_j M^b_{1/2}$, then from $\hat{a}\hat{b}\hat{a}\hat{c}\hat{b}$ there is an impossible square in $L$ on $M^a_1 \cdots M^a_j M^b_{1/2} \cdots M^b_k$.  Therefore, this $\hat{c}$ spans a block in
\[
M^a_1 \cdots M^a_j M^{c,1}_{1/2} \cdots M^{c,1}_l M^b_{1/2} \cdots M^b_k.
\]
That means $|\hat{c}| \geq 13$ and $\hat{c}$ is rigid.  From $\hat{c}\hat{a}$ we see that $\hat{c}$ is right flush, which means $\hat{b}$ is blush.  From $\hat{b}\hat{c}$ we see that $\hat{c}$ is blush,  and from $\hat{a}\hat{b}$, so is the locally rigid $\hat{a}$, which now is also seen to be rigid with length a multiple of 13.  This completes the proof starting with $\hat{a}$ in the principal role.

\medskip

Next assume that $\hat{b}$ is locally rigid in this local instance of $\kappa_2$ and left flush.  If $\hat{a}\hat{c}$ is contained within $M^b_k M^b_1$, then there is an impossible square on $M^b_{1/2} \cdots M^b_k$.  Therefore, $\hat{a}\hat{c}$ is at least a block and has length $\geq 13$.  That means that either $|\hat{a}| \geq 7$ or $|\hat{c}| \geq 7$.  If $|\hat{a}| = 7$, then $\hat{a}$ is locally rigid, because Theorem 1 offers no hope if it is not (and besides, $\hat{b}$ is locally rigid).  If $|\hat{c}| = 7$,  then $\hat{c}$ is locally rigid, because Theorem 1 is no shelter, with $\hat{a} \geq 6$ (and besides, $\hat{b}$ is locally rigid). We conclude at least one of $\hat{a}$ or $\hat{c}$ is locally rigid and maybe more.

We start with $\hat{a}$ locally rigid and right flush (from $\hat{a}\hat{b}$).  So if the left $\hat{c}$ goes no further than $M^b_{1/2}$, there is an impossible square on $M^a_1 \cdots M^a_j M^b_{1/2} \cdots M^b_k$ in $L$ from $\hat{a}\hat{b}\hat{a}\hat{c}\hat{b}$.  Therefore, this $\hat{c}$ spans a block, has length at least 13, and is rigid and blush (from $\hat{a}\hat{c}\hat{b})$.  So the locally rigid $\hat{a}$ is blush from $\hat{c}\hat{a}$ and indeed rigid, and $\hat{b}$ is blush from $\hat{b}\hat{a}$ and hence rigid.

We turn to $\hat{c}$ locally rigid and right flush (from $\hat{c}\hat{b}$).  If $\hat{a}$ fits in $M^c_1$, then from $\hat{b}\hat{a}\hat{c}\hat{b}\hat{c}$, $L$ has the impossible square on $M^b_1\cdots M^b_k M^c_{1/2} \cdots M^c_l$.  Thus $|\hat{a}| \geq 13$, and $\hat{a}$ is rigid and blush from $\hat{c}\hat{a}\hat{b}$.  So $\hat{c}$ is blush from $\hat{a}\hat{c}$ and rigid from crossing a block or more, and the locally rigid $\hat{b}$ is blush from $\hat{b}\hat{c}$ and rigid.  This completes the proof starting with $\hat{b}$ in the principal role.

\medskip

Let's assume that $\hat{c}$ is locally rigid in this local instance of $\kappa_2$ and left flush. The argument follows the proof for $\hat{a}$ above.

\medskip

Finally, we note that the cases for the keyword being right flush are exactly symmetric by palindromicity.   \qed

\bigskip

{\bf Theorem 4  (Shorter, nonrigid instance of $\kappa_2$).}  {\it  Suppose $[\hat{a},\hat{b},\hat{c},\rho]$ is a local instance of $\kappa_2$ with a locally rigid keyword.  Then there is a shorter local instance $[\dot{a},\dot{b},\dot{c}, \rho']$ with no rigidity and $|\rho'| \leq \lfloor (|\rho|+6)/13\rfloor$.} 

\medskip
{\it Proof.}  Let $[\hat{a},\hat{b},\hat{c},\rho]$ be a local instance of $\kappa_2$ with a locally rigid keyword and length $k$.  By the two previous theorems, we have an assignment of the same length for a local instance $[ \check{a},\check{b},\check{c}, \tau]$ where all keywords are sequences of blocks, and $|\tau| \leq|\rho|+6$.  Let $\mu =  \langle \check{a},\check{b},\check{c} \rangle$; $\mu$ is its own matrix and $|\mu| = k$.  

Let $\nu = \rho'\mu'$ be an initial subterm of $L$ such that $\overline{\mu'} = \mu$ and $\overline{\rho'}= \tau$; thus $|\nu| = (|\tau|+k)/13$ and $|\rho'| \leq \lfloor (|\rho|+6)/13\rfloor$.  Then $\mu'$ determines keywords for a local instance $[\dot{a},\dot{b},\dot{c}, \rho' ]$ such that $\overline{\dot{a}} = \check{a}, \overline{\dot{b}} = \check{b},\overline{\dot{c}} = \check{c}$ and $|\langle \dot{a},\dot{b},\dot{c} \rangle| = k/13$.  If the new local instance has local rigidity, repeat the procedure until it terminates, which will happen because the newer instance is strictly shorter than the older.    \qed

\bigskip

To estimate the number of rounds in Theorem 4, calculate $\lfloor \log_{13} k\rfloor$.  When there is  no rigidity, all keywords are length 7 or less, and the length of the instance is no more than $9\times 6$; in the only cases when a 7-letter keyword can appear, the length is 21.  
\bigskip

{\bf \S 3.  Results on short keywords.}

\medskip

We turn to the analysis of keyword candidates that are not locally rigid.  That immediately limits attention to 7 letters or fewer, and we have already taken care of 7.  In each case, we assume a local instance $[ \hat{a},\hat{b},\hat{c},\rho]$ in $L$ where no keyword is locally rigid.  For all keywords of 6 or fewer letters, we find a contradiction.  We start with 6-letter, 5-letter, and 4-letter.  Finally, we quickly handle the shortest feasible keywords.

Among the fourteen 6-letter subwords of $L$, up to rot variants, the ten rigid ones are:  
$ABACAB, ABACBA, ABACBC, ABCABA, ABCBAB,$ $ACABAC,$ $ACABCA,$ $ACABCB, ACBACA, ACBCAC$

\bigskip

{\bf Theorem 5 (Six-letter keywords).}  {\it There is no local instance $[\hat{a},\hat{b},\hat{c},\rho]$ of $\kappa_2$ such that no keyword is locally rigid and the longest keyword has 6 letters.}

\medskip

{\it Proof.} The words we need to cover are listed, followed by their matrices.

\[
\begin{array}{c}
ABCACB\ \  (\overline{AB}, \overline{B})\\
ABCBAC\ \ (\overline{AC}, \overline{A})\\ 
ACBABC\ \ (\overline{C}, \overline{BA})\\
ACBCAB\ \ (\overline{BC}, \overline{A})\\ 
\end{array}
\]

\medskip

{\bf ABCACB}\ \  $(\overline{AB}, \overline{B})$

(I)  We investigate $\langle \hat{a}, \hat{b}, \hat{c} \rangle = \langle ABCACB, -, - \rangle$ where two adjacent local occurrences of $\hat{a} = ABCACB$ are not rigid (because $\hat{a}$ is not locally rigid, such a pair exists).  

(A)  The first case is  $ABCACB x ABCACB$ with matrix $\overline{AB} \cdots \overline{B}$ (the shortest example is when there is just $\overline{AB}$ and $x = AC$). 

(1) So if this is an instance of $\hat{a} \hat{b} \hat{a}$, then $x = \hat{b}$ has the form $AC$, because the alternative $ACABC\cdots ACBAC$ would be rigid.  The entire string $\hat{a} \hat{b} \hat{a}$ is rigid, because its length is at least 14.  Thus, whether we started right or left with $\hat{a} \hat{b} \hat{a}$, $\hat{c}$ is left flush in $\hat{a}\hat{b} \hat{a}\hat{c}$, and therefore $\hat{c}$ begins like $\overline{A} = ABCB\cdots$ to avoid end-letter squares; note that $\hat{c} = AB$ or $ABC$ is not adequate.  So $\hat{b}\hat{c} = ACABCB\cdots$ (rigid), and $\hat{c}$ will need at least 12 letters to get to the next $\hat{a}$ in $\hat{b}\hat{c}\hat{a}$, making $\hat{c}$ rigid.

(2)  If this is an instance of $\hat{a}\hat{c}\hat{b}\hat{c}\hat{a}$, then $x = \hat{c}\hat{b}\hat{c} = ACABCA\cdots CACBAC$, at least 28 letters.

(a)  If $\hat{c} = AC$, then $\hat{b}$ has at least 24 letters and is rigid.

(b)  If $\hat{c} \neq AC$, then $\hat{c}$ has at least 10 letters and is rigid.

(B)  The second case is $ABCACB x ABCACB$ with matrix $\overline{B} \cdots \overline{AB}$ (the shortest example is when there is just $\overline{BAB}$ and $x = ABCBA\ CBC\ ABCB$). 

(1)  In the case $x = \hat{b}$, $|\hat{b}| \geq 12$, so it is rigid.

(2)  Consider the case $x = \hat{c}\hat{b}\hat{c}$.

(a)  If $\hat{c} = ABCB$, then $\hat{b} = ACBC\cdots$, so $\hat{a}\hat{b}$ has a square on $ACB$.

(b)   Otherwise, $\hat{c}$ has at least 12 letters and is rigid.

\medskip

(II)  Next we study $\langle \hat{a}, \hat{b}, \hat{c} \rangle = \langle -,ABCACB,  - \rangle$ where two adjacent local occurrences of $\hat{b} = ABCACB$ are not rigid (if $\hat{b}$ is not locally rigid, then such a pair exists).  

(A)  First we examine  $ABCACB x ABCACB$ with matrix $\overline{AB} \cdots \overline{B}$ (the shortest example is when there is just $\overline{AB}$ and $x = AC$).

(1) If $x = AC$, then one of $\hat{a}$ and $\hat{c}$ is $A$ and the other is $C$.  This leads to a square on $A$.

(2)  If $x = ACAB\cdots CBAC$, no matter how it is split, there will be a square on $AC$ (and $CA$)when $\hat{a}$ and $\hat{c}$ are reversed.

(B)  The second case is $ABCACB x ABCACB$ with matrix $\overline{B} \cdots \overline{AB} $ (the shortest example is when there is just $\overline{BAB}$ and $x = ABCBA\,CBC\,ABCB$).  One of $\hat{a}$ or $\hat{c}$ has disallowed length $\geq 7$, because the even 6-letter split of the shortest $x$ into $ABCBAC$ and $BCABCB$ gives a square on $B$ with $\hat{b}$.

\medskip

(III)   Finally we check $\langle \hat{a}, \hat{b}, \hat{c} \rangle = \langle -,-,ABCACB \rangle$.  

(A)  The first case is  $ABCACB x ABCACB$ with matrix $\overline{AB}\cdots \overline{B}$ (the shortest example is when there is just $\overline{AB}$ and $x = AC$).  We note $x = \hat{b}$ has the form $AC$ only.  Furthermore, $\hat{a}$ must be left flush in $\hat{b}\hat{a}\hat{c}$).   Therefore $\hat{a}$ begins like $\overline{A} = ABCB\cdots$ to avoid end-letter squares.  So $\hat{b}\hat{a} = ACABCB\cdots$ (rigid), and $\hat{a}$ will need at least 12 letters to get to the next $\hat{c}$, making it rigid.  

(B)  The second case is $ABCACB x ABCACB$ with matrix $\overline{B} \cdots \overline{AB}$ (the shortest example is when there is just $\overline{BAB}$ and $x = ABCBA\,CBC\,ABCB$).  Thus, $x = \hat{b}$, $|\hat{b}| \geq 12$, so it is rigid.

This completes the proof for the 6-letter nonrigid $ABCACB$.

\bigskip

{\bf ABCBAC}\ \  $(\overline{A}, \overline{AC})$

\medskip

(I)  We investigate $\langle \hat{a}, \hat{b}, \hat{c} \rangle = \langle ABCBAC, -, - \rangle$ where two adjacent local occurrences of $\hat{a} = ABCBAC$ are not rigid (because $\hat{a}$ is not locally rigid, such a pair exists).  

(A)  The first case is  $ABCBAC x ABCBAC$ with matrix $\overline{A} \cdots \overline{AC}$ (the shortest example is when there is just $\overline{AC}$ and $x = BC$). 

(1) So if this is an instance of $\hat{a} \hat{b} \hat{a}$, then $x = \hat{b}$ has the form $BC$ for the form $BCABCB\cdots CBACBC$ is too long.  The entire string $\hat{a} \hat{b} \hat{a}$ is rigid, because its length is at least 14.  Thus, $\hat{c}$ is right flush in $\hat{c} \hat{a} \hat{b} \hat{a}$ and therefore $\hat{c}$ ends like $\overline{C} = \cdots ABAC$ to avoid end-letter squares; note that $\hat{c} = AC$ or $BAC$ is not adequate.  So $\hat{c}\hat{b} = \cdots ABACBC\cdots$ (rigid), and $\hat{c}$ will need at least 12 letters to go left to the next $\hat{a}$ in $\hat{b}\hat{a}\hat{c}$, making $\hat{c}$ rigid.

(2)  If this is an instance of $\hat{a}\hat{c}\hat{b}\hat{c}\hat{a}$, then $x = \hat{c}\hat{b}\hat{c} = BCABCB\cdots CBACBC$.

(a)  If $\hat{c} = BC$, then $\hat{b}$ has at least 24 letters and is rigid.

(b)  If $\hat{c} \neq BC$, then $\hat{c}$ has at least 10 letters and is rigid.

(B)  The second case is $ABCBAC x ABCBAC$ with matrix $\overline{AC} \cdots \overline{A}$ (the shortest example is when there is just $\overline{ACA}$ and $x = ABAC\ BAB\ CABAC$). 

(1)  In the case $x = \hat{b}$, $|\hat{b}| \geq 12$, so it is rigid.

(2)  Consider the case $x = \hat{c}\hat{b}\hat{c}$.

(a)  If $\hat{c} = ABAC$, then $\hat{b} = BABC\cdots$, which does not fit after either instance of $\hat{a}$.

(b)   Otherwise, $\hat{c}$ has at least 12 letters and is rigid.

\medskip

(II)  Next we study $\langle \hat{a}, \hat{b}, \hat{c} \rangle = \langle -,ABCBAC,  - \rangle$ where two adjacent local occurrences of $\hat{b} = ABCBAC$ are not rigid ($\hat{b}$ is not locally rigid, so such a pair exists).  

(A)  First we examine   $ABCBAC x ABCBAC$ with matrix $\overline{A} \cdots \overline{AC}$ (the shortest example is when there is just $\overline{AC}$ and $x = BC$). 

(1) If $x = BC$, then one of $\hat{a}$ and $\hat{c}$ is $B$ and the other is $C$.  This leads to a square on $A$.

(2)  If $x = BCAB\cdots ACBC$, no matter how it is split, there will be a square on $BC$ (and $CB$) when $\hat{a}$ and $\hat{c}$ are reversed.

(B)  The second case is $ABCBAC x ABCBAC$ with matrix $\overline{AC}\cdots \overline{A}$ (the shortest example is when there is just $\overline{ACA}$ and $x = ABAC\ BAB\ CABAC$).  One of $\hat{a}$ or $\hat{c}$ has illegal length $\geq 7$, because the even split of the shortest $x$ into $ABACBA$ and $BCABAC$ gives a square on $B$ with $\hat{b}$.

\medskip

(III)   Finally we check $\langle \hat{a}, \hat{b}, \hat{c} \rangle = \langle -,-,ABCBAC \rangle$.

(A)  The first case is  $ABCBAC x ABCBAC$ with matrix $\overline{A} \cdots \overline{AC}$ (the shortest example is when there is just $\overline{AC}$ and $x = BC$).  We note that $x = \hat{b}$ has the form $BC$, the alternative being too long.  Furthermore, $\hat{a}$ must be right flush.  Therefore $\hat{a}$ ends like $\overline{C} = \cdots ABAC$ to avoid end-letter squares.  So $\hat{a}\hat{b} = \cdots ABACBC\cdots$ (rigid), and $\hat{a}$ will need at least 12 letters to go left to the next $\hat{c}$ in $\hat{c}\hat{a}\hat{b}$, making $\hat{a}$ rigid.

(B)   The second case is $ABCBAC x ABCBAC$ with matrix $\overline{AC} \cdots \overline{A}$ (the shortest example is when there is just $\overline{ACA}$ and $x = ABAC\ BAB\ CABAC$).   Thus, $x = \hat{b}$ and has at least 12 letters, so it is rigid.

This completes the proof for the 6-letter nonrigid $ABCBAC$.

\bigskip

{\bf ACBABC}\ \ $(\overline{BA}, \overline{C})$

(I)  We investigate $\langle \hat{a}, \hat{b}, \hat{c} \rangle = \langle ACBABC, -, - \rangle$ where two adjacent local occurrences of $\hat{a} = ACBABC$ are not rigid ($\hat{a}$ is not locally rigid, so such a pair exists).  

(A)  The first case is  $ACBABC x ACBABC$ with matrix $\overline{BA} \cdots \overline{C}$ (the shortest example is when there is just $\overline{BAC}$ and $x = CBA\ CBC\ ABCBA\ CAB$). 

(1) So if this is an instance of $\hat{a} \hat{b} \hat{a}$, then $|x| = |\hat{b}| \geq 14$ means $\hat{b}$ is rigid.

(2)  If this is an instance of $\hat{a}\hat{c}\hat{b}\hat{c}\hat{a}$, then $x = \hat{c}\hat{b}\hat{c} = CBAC\cdots ACAB$.  This means $|\hat{c}| \geq 8$, and $\hat{c}$ is rigid.

(B)  The second case is $ACBABC x ACBABC$ with matrix $\overline{C} \cdots \overline{BA}$ (the shortest example is when there is just $\overline{CBA}$ and $x = ABAC\ BCACB\ ACA\ BC$). 

(1)  In the case $x = \hat{b}$, $|\hat{b}| \geq 14$, so it is rigid.

(2)  In the case $x = \hat{c}\hat{b}\hat{c}$, $|\hat{c}| \geq 8$, and $\hat{c}$ is rigid.

\medskip

(II)  Next we study $\langle \hat{a}, \hat{b}, \hat{c} \rangle = \langle -,ACBABC,  - \rangle$ where two adjacent local occurrences of $\hat{b} = ACBABC$ are not rigid.  

(A)  First we examine   $ACBABC x ACBABC$ with matrix $\overline{BA} \cdots \overline{C}$ (the shortest example is when there is just $\overline{BAC}$ and $x = CBA\ CBC\ ABCBA\ CAB$).  One of $\hat{a}$ and $\hat{c}$ has at least 7 letters, which contradicts hypotheses.  

(B)  The second case is $ACBABC x ACBABC$ with matrix $\overline{C}  \cdots \overline{BA}$ (the shortest example is when there is just $\overline{CBC}$ and $x = ABAC\ BCACB\ ACA\ BC$).  One of $\hat{a}$ or $\hat{c}$ has a disallowed length $\geq 7$.

\medskip

(III)   Finally we check $\langle \hat{a}, \hat{b}, \hat{c} \rangle = \langle -,-,ACBABC \rangle$. 

(A)  The first case is  $ACBABC x ACBABC$ with matrix $\overline{BA} \cdots \overline{C}$ (the shortest example is when there is just $\overline{BAC}$ and $x = CBA\ CBC\ ABCBA\ CAB$).    
Then $|x| = |\hat{b}| \geq 14$ means $\hat{b}$ is rigid.

(B)  The second case is $ACBABC x ACBABC$ with matrix $\overline{C}  \cdots \overline{BA}$ (the shortest example is when there is just $\overline{CBA}$ and $x = ABAC\ BCACB\ ACA\ BC$).  Again, $\hat{b}$ has at least 14, so it is rigid.

This completes the proof for the 6-letter nonrigid $ACBABC$.

\bigskip

{\bf ACBCAB}\ \ $(\overline{A}, \overline{BC})$

(I)  We investigate $\langle \hat{a}, \hat{b}, \hat{c} \rangle = \langle ACBCAB, -, - \rangle$ where two adjacent local occurrences of $\hat{a} = ACBCAB$ are not rigid ($\hat{a}$ is not locally rigid, so such a pair exists).  

(A)  The first case is  $ACBCAB x ACBCAB$ with matrix $\overline{A} \cdots \overline{BC}$ (the shortest example is when there is just $\overline{ABC}$ and $x = CBA\ BCACB\ ACA\ BC$). 

(1) So if this is an instance of $\hat{a} \hat{b} \hat{a}$, then $|x| = |\hat{b}| \geq 13$ means $\hat{b}$ is rigid.

(2)  If this is an instance of $\hat{a}\hat{c}\hat{b}\hat{c}\hat{a}$, then $x = \hat{c}\hat{b}\hat{c} = CBAB\cdots CABC$.  

(a)  If $\hat{c} = C$, then $|\hat{b} \geq 11$, and $\hat{b}$ is rigid.

(b)  Otherwise, $|\hat{c}| \geq 8$, and $\hat{c}$ is rigid.

(B)  The second case is $ACBCAB x ACBCAB$ with matrix $\overline{BC}\cdots \overline{A}$ (the shortest example is when there is just $\overline{CBA}$ and $x = AC\ BAB\ CABAC\ ABCB$). 

(1)  In the case $x = \hat{b}$, $|\hat{b}| \geq 14$, so $\hat{b}$ is rigid.

(2)  In the case $x = \hat{c}\hat{b}\hat{c}$, $|\hat{c}| \geq 7$, not allowed.

\medskip

(II)  Next we study $\langle \hat{a}, \hat{b}, \hat{c} \rangle = \langle -,ACBCAB,  - \rangle$ where two adjacent local occurrences of $\hat{b} = ACBCAB$ are not rigid.  

(A)  First we examine   $ACBCAB x ACBCAB$ with matrix $\overline{A} \cdots\overline{BC}$ (the shortest example is when there is just $\overline{ABC}$ and $x = CBA\ BCACB\ ACA\ BC$).  One of $\hat{a}$ and $\hat{c}$ has at least 7 letters, which is disallowed.  

(B)  The second case is $ACBABC x ACBABC$ with matrix $\overline{BC} \cdots \overline{A}$ (the shortest example is when there is just $\overline{BCA}$ and $x = AC\ BAB\ CABAC\ ABCB$).  One of $\hat{a}$ and $\hat{c}$ has length $\geq 7$ and that is disallowed.

\medskip

(III)   Finally we check $\langle \hat{a}, \hat{b}, \hat{c} \rangle = \langle -,-,ACBCAB \rangle$. 

(A)  The first case is  $ACBCAB x ACBCAB$ with matrix $\overline{A} \cdots\overline{BC}$ (the shortest example is when there is just $\overline{ABC}$ and $x = CBA\ BCACB\ ACA\ BC$).     
Then $|x| = |\hat{b}| \geq 13$ means $\hat{b}$ is rigid.

(B)  The second case is $ACBCAB x ACBCAB$ with matrix $\overline{BC} \cdots \overline{A}$ (the shortest example is when there is just $\overline{BCA}$ and $x = AC\ BAB\ CABAC\ ABCB$).  Thus, $x = \hat{b}$, $|\hat{b}| \geq 14$, so it is rigid.

This completes the proof for the 6-letter nonrigid $ACBCAB$.

We conclude that no nonrigid instance of $\kappa_2$ in $L$ has a longest keyword with six (or more) letters.   \qed

\bigskip

Of the ten 5-letter subwords of $L$ up to rotation, six are not rigid, and four are:  $ABACA$, $ABCAB$, $ACABA$, $ACBAC$.

\medskip

{\bf Theorem 6 (Five-letter keywords).}  {\it There is no local instance $[\hat{a},\hat{b},\hat{c},\rho]$ of $\kappa_2$ such that no keyword is locally rigid and the longest keyword has 5 letters.}

\medskip

{\it Proof.} The six nonrigid 5-letter words we need to cover are listed, followed by their matrices.

\[
\begin{array}{ll}
ABACB &  (\overline{C}, \overline{CB})\\
ABCAC &(\overline{AB}, \overline{A})\\ 
ABCBA & (\overline{A}, \overline{A})\\
ACABC & (\overline{B}, \overline{CA})\\
ACBAB & (\overline{C}, \overline{BA})\\ 
ACBCA & (\overline{A},\overline{BC}, \overline{CB})\\ 
\end{array}
\]

\medskip

{\bf ABACB}\ \ $(\overline{C}, \overline{CB})$

(I)  We investigate $\langle \hat{a}, \hat{b}, \hat{c} \rangle = \langle ABACB, -, - \rangle$ where two adjacent local occurrences of $\hat{a} = ABACB$ are not rigid (because $\hat{a}$ is not locally rigid, such a pair exists). 
 
(A)  The first case is  $ABACB x ABACB$ with matrix $\overline{C} \cdots \overline{CB}$ (the shortest example is when there is just $\overline{CB}$ and $x = ABC$). 

(1) So if this is an instance of $\hat{a} \hat{b} \hat{a}$, then $x = \hat{b} = ABC$ or $ABCAB \cdots CBABC$.

(a)  If $\hat{b} = ABC$, then it must be that $\hat{c} = AC$ or $AC \cdots  AC$ to avoid squares (in particular, if $\hat{c} = AB\cdots C$, then there is a square on $CAB$ in $\hat{c}\hat{b}\hat{c}$).

(i)  If $\hat{c} = AC$, then there is a square on $CA$ in $\hat{a}\hat{c}\hat{b}$.

(ii)  If $\hat{c} = AC \cdots AC$, then the third letter cannot be $A$ because that gives a square on $CA$ in $\hat{b}\hat{c}$, and it cannot be $B$ because that gives a square on $ACB$ in $\hat{a}\hat{c}$.

(b)  If $\hat{b} = ABCAB \cdots CBABC$, then $|\hat{b}| \geq 8$, and $\hat{b}$ is rigid.

(2)  If this is an instance of $\hat{a}\hat{c}\hat{b}\hat{c}\hat{a}$, then $x = \hat{c}\hat{b}\hat{c} = ABCAB \cdots CBABC$.  $ABC$ is not an option.  Therefore, $|\hat{c}| \geq 8$, and $\hat{c}$ is rigid.

(B)  The second case is $ABACB x ABACB$ with matrix $\overline{CB}\cdots \overline{C}$ (the shortest example is when there is just $\overline{CBC}$ and $x = CACB\ ACA\ BCACB\ C$). 

(1)  In the case $x = \hat{b}$, $|\hat{b}| \geq 13$, so it is rigid.

(2)  In the case $x = \hat{c}\hat{b}\hat{c}$, $\hat{c} = C$ or $CACBA \cdots CACBC$.

(a)  If $\hat{c} = C$, then $\hat{b} = ACBA \cdots ACBC$ and is rigid by length.

(b)  If $\hat{c} = CACBA \cdots CACBC$, then it is rigid by length.

\medskip

(II)  Next we study $\langle \hat{a}, \hat{b}, \hat{c} \rangle = \langle -,ABACB,  - \rangle$ where two adjacent local occurrences of $\hat{b} = ABACB$ are not rigid.  

(A)  First we examine  $ABACB x ABACB$ with matrix $\overline{C}\cdots  \overline{CB}$ (the shortest example is when there is just $\overline{CB}$ and $x = ABC$).   So $x = \hat{c}\hat{a}$ or $\hat{a}\hat{c}$ and has the form $ABC$ or $ABC \cdots ABC$.

(1)  If one of $\hat{a}, \hat{c}$ ends in $A$ or $AB$, there is a square with the beginning of $\hat{b}$.

(2)  Otherwise, one of $\hat{a}, \hat{c}$ starts with $ABC$ and the other ends with it.  That means when reversed there is a square on $ABC$.

(B)  The second case is $ABACB x ABACB$ with matrix $\overline{CB} \cdots \overline{C}$ (the shortest example is when there is just $\overline{CBC}$ and $x = CACB\ ACA\ BCACB\ C$).  No matter how $\hat{a}$ and $\hat{c}$ divide $x$, one starts with $C$ and the other ends with $C$, which gives a square on $C$ when they are reversed.

\medskip

(III)   Finally we check $\langle \hat{a}, \hat{b}, \hat{c} \rangle = \langle -,-,ABACB \rangle$. 

(A)  The first case is  $ABACB x ABACB$ with matrix $\overline{C}\cdots  \overline{CB}$ (the shortest example is when there is just $\overline{CB}$ and $x = ABC$).  So $x = \hat{b} = ABC$ or $ABCAB \cdots CBABC$. 

(1)  If $\hat{b} = ABC$, then it must be that $\hat{a} = AC$ or $AC \cdots  AC$ to avoid squares.

(a)  If $\hat{a} = AC$, then there is a square on $CA$ in $\hat{b}\hat{a}\hat{c}$.

(b)  If $\hat{a} = AC \cdots AC$, then the third letter cannot be $A$ because that gives a square on $CA$ in $\hat{b}\hat{a}$, and it cannot be $B$ because that gives a square on $ACB$ in $\hat{c}\hat{a}$.

(2)  If $\hat{b} = ABCAB \cdots CBABC$, then $|\hat{b}| \geq 8$, and $\hat{b}$ is rigid.

(B)  The second case is $ABACB x ABACB$ with matrix $\overline{CB}\cdots \overline{C}$ (the shortest example is when there is just $\overline{CBC}$ and $x = CACB\ ACA\ BCACB\ C$).  Because $x = \hat{b}$ has at least 13 letters, it is rigid.

\bigskip

{\bf ABCAC}\ \ $(\overline{AB}, \overline{B})$

(I)  We investigate $\langle \hat{a}, \hat{b}, \hat{c} \rangle = \langle ABCAC, -, - \rangle$ where two adjacent local occurrences of $\hat{a} = ABCAC$ are not rigid (since $\hat{a}$ is not locally rigid, such a pair exists). 
 
(A)  The first case is  $ABCAC x ABCAC$ with matrix $\overline{AB} \cdots \overline{B}$ (the shortest example is when there is just $\overline{AB}$ and $x = BAC$). 

(1) So if this is an instance of $\hat{a} \hat{b} \hat{a}$, then $x = \hat{b} = BAC$ or $BACAB \cdots ACBAC$.

(a)  If $\hat{b} = BAC$, then it must be that $\hat{c} = BC$ or $BC \cdots  BC$ to avoid squares (in particular, if $\hat{c} = BA\cdots C$, then there is a square on $CBA$ in $\hat{c}\hat{b}\hat{c}$).

(i)  If $\hat{c} = BC$, then there is a square on $CB$ in $\hat{a}\hat{c}\hat{b}$.

(ii)  If $\hat{c} = BC \cdots BC$, then the third letter from the end cannot be $A$ because that gives a square on $ABC$ in $\hat{c}\hat{a}$, and it cannot be $C$ because that gives a square on $CB$ in $\hat{c}\hat{b}$.

(b)  If $\hat{b} = BACAB \cdots ACBAC$, then $|\hat{b}| \geq 8$, and $\hat{b}$ is rigid.

(2)  If this is an instance of $\hat{a}\hat{c}\hat{b}\hat{c}\hat{a}$, then $x = \hat{c}\hat{b}\hat{c} =  BACAB \cdots ACBAC$.  $BAC$ is not an option.  Therefore, $|\hat{c}| \geq 8$, and $\hat{c}$ is rigid.

(B)  The second case is $ABCAC x ABCAC$ with matrix $\overline{B} \cdots \overline{AB} $ (the shortest example is when there is just $\overline{BAB}$ and $x = B\ ABCBA\ CBC\ ABCB$). 

(1)  In the case $x = \hat{b}$, $|\hat{b}| \geq 13$, so it is rigid.

(2)  In the case $x = \hat{c}\hat{b}\hat{c}$, $\hat{c} = B$ or $BABCB \cdots CABCB$.

(a)  If $\hat{c} = B$, then $\hat{b} = ABCB \cdots CABC$ and is rigid by length.

(b)  If $\hat{c} = BABCB \cdots CABCB$, then it is rigid by length.

\medskip

(II)  Next we study $\langle \hat{a}, \hat{b}, \hat{c} \rangle = \langle -,ABCAC,  - \rangle$ where two adjacent local occurrences of $\hat{b} = ABCAC$ are not rigid.  

(A)  First we examine  $ABCAC x ABCAC$ with matrix $\overline{AB} \cdots \overline{B}$ (the shortest example is when there is just $\overline{AB}$ and $x = BAC$).  So $x = \hat{c}\hat{a}$ or $\hat{a}\hat{c}$ and has the form $BAC$ or $BAC \cdots BAC$.

(1)  If one of $\hat{a}, \hat{c}$ ends in $A$ or $BA$, there is a square with the start of $\hat{b}$.

(2)  Otherwise, one of $\hat{a}, \hat{c}$ starts with $BAC$ and the other ends with it.  That means when reversed they give a square on $BAC$.

(B)  The second case is $ABCAC x ABCAC$ with matrix $\overline{B} \cdots \overline{AB}$ (the shortest example is when there is just $\overline{BAB}$ and $x = B\ ABCBA\ CBC\ ABCB$).  No matter how $\hat{a}$ and $\hat{c}$ divide $x$, one starts with $B$ and the other ends with $B$, which gives a square on $B$ when they are reversed.

\medskip

(III)   Finally we check $\langle \hat{a}, \hat{b}, \hat{c} \rangle = \langle -,-,ABCAC \rangle$. 

(A)  The first case is  $ABCAC x ABCAC$ with matrix $\overline{AB} \cdots \overline{B}$ (the shortest example is when there is just $\overline{AB}$ and $x = BAC$).   So $x = \hat{b} = BAC$ or $BACA \cdots CBAC$. 

(1)  If $\hat{b} = BAC$, then it must be that $\hat{a} = BC$ or $BC \cdots  BC$ to avoid squares (in particular, if $\hat{a} = BA\cdots C$, then there is a square on $CBA$ in $\hat{a}\hat{b}\hat{a}$).

(a)  If $\hat{a} = BC$, then there is a square on $CB$ in $\hat{c}\hat{a}\hat{b}$.

(b)  If $\hat{a} = BC \cdots BC$, then the third letter from the end cannot be $A$ because that gives a square on $ABC$ in $\hat{a}\hat{c}$, and it cannot be $C$ because that gives a square on $CB$ in $\hat{a}\hat{b}$.

(2)  If $\hat{b} = BACA \cdots CBAC$, then $|\hat{b}| \geq 8$, and $\hat{b}$ is rigid.

(B)  The second case is $ABCAC x ABCAC$ with matrix $\overline{B}\cdots \overline{AB} $ (the shortest example is when there is just $\overline{BAB}$ and $x = B\ ABCBA\ CBC\ ABCB$).  Because $x = \hat{b}$ has at least 13 letters, it is rigid.

\bigskip

{\bf ABCBA}\ \ $(\overline{A}, \overline{A})$

(I)  We investigate $\langle \hat{a}, \hat{b}, \hat{c} \rangle = \langle ABCBA, -, - \rangle$ where two adjacent local occurrences of $\hat{a} = ABCBA$ are not rigid. 
 
(A)  The first case is  $ABCBA x ABCBA$ with matrix $\overline{A} \cdots \overline{A}$ (the shortest example is when there is just $\overline{A}$ and $x = CBC$). 

(1) So if this is an instance of $\hat{a} \hat{b} \hat{a}$, then $x = \hat{b} = CBC$ or $CBC \cdots CBC$. So $\hat{c}$ must begin and end with $B$.  That means $\hat{b}\hat{c}$ has a square on $BC$ (and $\hat{c}\hat{b}$ has a square on $CB$).

(2)  If this is an instance of $\hat{a}\hat{c}\hat{b}\hat{c}\hat{a}$, then $x = \hat{c}\hat{b}\hat{c} = CBC$ or $CBCAB \cdots BACBC$ ($CBCACBC, CBCABACBC$ are not subwords of $L$). 

(a)  If $\hat{c} = C$ and $\hat{b} = B$, then $\hat{a}\hat{b}\hat{a}$ includes the subword $BABAB$.

(b)  If $\hat{c} = C$ and $\hat{b} = BCAB\cdots BACB$, then $\hat{b}$ has 8 or more letters and is rigid.

(c)  If $\hat{c} = CBC$, then $\hat{b}$ starts and ends with $A$, square to $\hat{a}$.

(B)  The second case is $ABCBA x ABCBA$ with matrix $\overline{A} \cdots \overline{A}$ where the two instances of $\hat{a}$ are right and left flush, respectively.  So $x$ is a string of one or more blocks.
  
(1)  If $x$ is $\hat{b}$, then $\hat{b}$ is rigid by length.

(2)  Suppose $x$ is $\hat{c}\hat{b}\hat{c}$.

(a)  If $x$ is more than one block, then either $\hat{b}$ or $\hat{c}$ is rigid by length---one of them has at least 9 letters.

(b)  Suppose $x$ is a single block.  It could be $\overline{C}$ (or it could be $\overline{B}$).  Note that because $\hat{c}$ is not rigid, we cannot conclude that it is blush.

(i)  If $\hat{c} = C$, then $\hat{b}$ has 11 letters and is rigid.

(ii)  If $\hat{c} = CABAC$, then $\hat{b} = BAB$, which is inconsistent with $\hat{b}$ starting with, ending with, or just being $CBC$ by lying next to $\hat{a}$.

\medskip

(II)  Next we study $\langle \hat{a}, \hat{b}, \hat{c} \rangle = \langle -,ABCBA,  - \rangle$ where two adjacent local occurrences of $\hat{b} = ABCBA$ are not rigid (if $\hat{b}$ is not locally rigid, then such a pair exists).  

(A)  First we examine  $ABCBA x ABCBA$ with matrix $\overline{A} \cdots \overline{A}$ (the shortest example is when there is just $\overline{A}$ and $x = CBC$).  So $x = CBC$ or $CBC \cdots CBC$. No matter how $x$ is divided into $\hat{a}$ and $\hat{c}$, one begins with $C$ and the other ends with it.  Reversing them gives a square on $C$.

(B)  The second case is $ABCBA x ABCBA$ with matrix $\overline{A} \cdots \overline{A}$ where the two instances of $\hat{b}$ are right and left flush. Whether $x$ is $\hat{a}\hat{c}$ or $\hat{c}\hat{a}$, we see that $x$ is a nonempty string of blocks, so one of $\hat{a}, \hat{c}$ has at least 7 letters, contrary to hypothesis.

\medskip

(III)   Finally we check $\langle \hat{a}, \hat{b}, \hat{c} \rangle = \langle -,-,ABCBA \rangle$. 

(A)  The first case is  $ABCBA x ABCBA$ with matrix $\overline{A} \cdots \overline{A}$ (the shortest example is when there is just $\overline{A}$ and $x = CBC$).  This is similar to (I.A.1) above.

(B)  The second case is $ABCBA x ABCBA$ with matrix $\overline{A} \cdots \overline{A}$ where the two instances of $\hat{c}$ are right and left flush.  So $x = \hat{b}$ consists of one or more blocks, and $\hat{b}$ is rigid by length.

\bigskip

{\bf ACABC}\ \ $(\overline{B}, \overline{CA})$

(I)  We investigate $\langle \hat{a}, \hat{b}, \hat{c} \rangle = \langle ACABC, -, - \rangle$ where two adjacent local occurrences of $\hat{a} = ACABC$ are not rigid.
 
(A)  The first case is  $ACABC x ACABC$ with matrix $\overline{B} \cdots \overline{CA}$ (the shortest example is when there is just $\overline{BCA}$ and $x = ACB\ CABAC\ BAB\ CAB$). 

(1) So if this is an instance of $\hat{a} \hat{b} \hat{a}$, then $x = \hat{b}$ has at least 14 letters and is rigid.

(2)  If this is an instance of $\hat{a}\hat{c}\hat{b}\hat{c}\hat{a}$, then $x = \hat{c}\hat{b}\hat{c} = ACBC \cdots BCAB$.   Therefore, $|\hat{c}| \geq 8$, and $\hat{c}$ is rigid.

(B)  The second case is $ACABC x ACABC$ with matrix $\overline{CA}\cdots \overline{B}$ (the shortest example is when there is just $\overline{CAB}$ and $x = BA\ CBC\ ABCBA\ BCACB$). 

(1)  In the case $x = \hat{b}$, $|\hat{b}| \geq 15$, so it is rigid.

(2)  In the case $x = \hat{c}\hat{b}\hat{c}$, $\hat{c} = B$ or $BACBC \cdots BCACB$.

(a)  If $\hat{c} = B$, then $\hat{b} = ACBC \cdots BCAC$ and is rigid by length.

(b)  If $\hat{c} = BACB \cdots CACB$, then it is rigid by length.

\medskip

(II)  Next we study $\langle \hat{a}, \hat{b}, \hat{c} \rangle = \langle -,ACABC,  - \rangle$ where two adjacent local occurrences of $\hat{b} = ACABC$ are not rigid.

(A)  First we examine  $ACABC x ACABC$ with matrix $\overline{B} \cdots \overline{CA}$ (the shortest example is when there is just $\overline{BCA}$ and $x = ACB\ CABAC\ BAB\ CAB$).    One of $\hat{a}, \hat{c}$ has at least 7 letters---disallowed.

(B)  The second case is $ACABC x ACABC$ with matrix $\overline{CA}\cdots \overline{B}$ (the shortest example is when there is just $\overline{CAB}$ and $x = BA\ CBC\ ABCBA\ BCACB$).  One of $\hat{a}, \hat{c}$ has at least 8 letters, so it is rigid. 

\medskip

(III)   Finally we check $\langle \hat{a}, \hat{b}, \hat{c} \rangle = \langle -,-,ACABC \rangle$. 

(A)  The first case is  $ACABC x ACABC$ with matrix $\overline{B} \cdots \overline{CA}$ (the shortest example is when there is just $\overline{BCA}$ and $x = ACB\ CABAC\ BAB\ CAB$).  So  $x = \hat{b}$ has at least 14 letters and is rigid.

(B)  The second case is $ACABC x ACABC$ with matrix $\overline{CA} \cdots \overline{B}$ (the shortest example is when there is just $\overline{CAB}$ and $x = BA\ CBC\ ABCBA\ BCACB$).  So $x = \hat{b}$ has at least 15 letters, so it is rigid. 

\bigskip

{\bf ACBAB}\ \ $(\overline{BA}, \overline{C})$

(I)  We investigate $\langle \hat{a}, \hat{b}, \hat{c} \rangle = \langle ACBAB, -, - \rangle$ where two adjacent local occurrences of $\hat{a} = ACBAB$ are not rigid.
 
(A)  The first case is  $ACBAB x ACBAB$ with matrix $\overline{BA} \cdots \overline{C}$ (the shortest example is when there is just $\overline{BAC}$ and $x = CBA\ CBC\ ABCBA\ CAB$).  

(1) So if this is an instance of $\hat{a} \hat{b} \hat{a}$, then $x = \hat{b}$ has at least 14 letters and is rigid.

(2)  If this is an instance of $\hat{a}\hat{c}\hat{b}\hat{c}\hat{a}$, then $x = \hat{c}\hat{b}\hat{c} =  CBACBC \cdots AB$.  Therefore, $|\hat{c}| \geq 8$, and $\hat{c}$ is rigid.

(B)  The second case is $ACBAB x ACBAB$ with matrix $\overline{C} \cdots \overline{BA} $ (the shortest example is when there is just $\overline{CBA}$ and $x = CABAC\ BCACB\ ACA\ BC$). 

(1)  In the case $x = \hat{b}$, $|\hat{b}| \geq 15$, so it is rigid.

(2)  In the case $x = \hat{c}\hat{b}\hat{c}$, $\hat{c} = C$ or $CABAC \cdots ACABC$.

(a)  If $\hat{c} = C$, then $\hat{b} = ABAC \cdots ACAB$ and is rigid by length.

(b)  If $\hat{c} = CABAC \cdots ACABC$, then it is rigid by length.

\medskip

(II)  Next we study $\langle \hat{a}, \hat{b}, \hat{c} \rangle = \langle -,ACBAB,  - \rangle$ where two adjacent local occurrences of $\hat{b} = ABCAC$ are not rigid (if $\hat{b}$ is not locally rigid, then such a pair exists).  

(A)  First we examine  $ACBAB x ACBAB$ with matrix $\overline{BA} \cdots \overline{C}$ (the shortest example is when there is just $\overline{BAC}$ and $x = CBA\ CBC\ ABCBA\ CAB$).  One of $\hat{a}, \hat{c}$ has at least 7 letters, which is disallowed.

(B)  The second case is $ACBAB x ACBAB$ with matrix $\overline{C} \cdots \overline{BA}$ (the shortest example is when there is just $\overline{CBA}$ and $x = CABAC\ BCACB\ ACA\ BC$).  One of $\hat{a}, \hat{c}$ has at least 8 letters, so it is rigid.

\medskip

(III)   Finally we check $\langle \hat{a}, \hat{b}, \hat{c} \rangle = \langle -,-,ACBAB \rangle$. 

(A)  The first case is  $ACBAB x ACBAB$ with matrix $\overline{BA}\cdots \overline{C}$ (the shortest example is when there is just $\overline{BAC}$ and $x = CBA\ CBC\ ABCBA\ CAB$).  $x = \hat{b}$ has at least 14 letters and is rigid.

(B)  The second case is $ACBAB x ACBAB$ with matrix $\overline{C} \cdots \overline{BA}$ (the shortest example is when there is just $\overline{CBA}$ and $x = CABAC\ BCACB\ ACA\ BC$).  $x = \hat{b}$ has at least 15 letters and is rigid.

\bigskip

{\bf ACBCA}\ \ $(\overline{A}, \overline{BC},\overline{CB})$

This is unusual because there are three possible matrices and six matrix pairs.

\medskip

(I)  We investigate $\langle \hat{a}, \hat{b}, \hat{c} \rangle = \langle ACBCA, -, - \rangle$ where two adjacent local occurrences of $\hat{a} = ACBCA$ are not rigid (because $\hat{a}$ is not locally rigid, such a pair exists). 
 
(A)  The first case is  $ACBCA x ACBCA$ with matrix $\overline{A} \cdots \overline{BC}$ (the shortest example is when there is just $\overline{ABC}$ and $x = BCBA\ BCACB\ ACA\ BC$). 

(1) So if this is an instance of $\hat{a} \hat{b} \hat{a}$, then $x = \hat{b}$ has at least 14 letters and is rigid.

(2)  If this is an instance of $\hat{a}\hat{c}\hat{b}\hat{c}\hat{a}$, then $\hat{c} = BC$ or $BCBA \cdots CABC$.

(a)  If $\hat{c} = BC$, then $\hat{b}$ has at least 10 letters and is rigid.

(b)  If $\hat{c} = BCBA \cdots CABC$, then it has at least 8 letters and is rigid.

(B)  The second case is $ACBCA x ACBCA$ with matrix $\overline{BC} \cdots \overline{A}$ (the shortest example is when there is just $\overline{BCA}$ and $x = BAC\ BAB\ CABAC\ ABCB$).

(1) So if this is an instance of $\hat{a} \hat{b} \hat{a}$, then $x = \hat{b}$ has at least 15 letters and is rigid.

(2)  If this is an instance of $\hat{a}\hat{c}\hat{b}\hat{c}\hat{a}$, then $\hat{c} = B$ or $BACB \cdots ABCB$.

(a)  If $\hat{c} = B$, then $\hat{b}$ has at least 13 letters and is rigid.

(b)  If $\hat{c} = BCBA \cdots CABC$, then it has at least 8 letters and is rigid.

(C)  The third case is  $ACBCA x ACBCA$ with matrix $\overline{A} \cdots \overline{CB}$ (the shortest example is when there is just $\overline{ACB}$ and $x = BCBA\ CABAC\ BAB\ CAB$). 

(1) So if this is an instance of $\hat{a} \hat{b} \hat{a}$, then $x = \hat{b}$ has at least 15 letters and is rigid.

(2)  If this is an instance of $\hat{a}\hat{c}\hat{b}\hat{c}\hat{a}$, then $\hat{c} = B$ or $BCBA \cdots BCAB$.

(a)  If $\hat{c} = B$, then $\hat{b}$ has at least 13 letters and is rigid.

(b)  If $\hat{c} = BCBA \cdots CABC$, then it has at least 8 letters and is rigid.

(D)  The fourth case is $ACBCA x ACBCA$ with matrix $\overline{CB} \cdots \overline{A}$ (the shortest example is when there is just $\overline{CBA}$ and $x = CB\ ACA\ BCACB\ ABCB$). 

(1) So if this is an instance of $\hat{a} \hat{b} \hat{a}$, then $x = \hat{b}$ has at least 14 letters and is rigid.

(2)  If this is an instance of $\hat{a}\hat{c}\hat{b}\hat{c}\hat{a}$, then $\hat{c} = CB$ or $CBAC \cdots ABCB$.

(a)  If $\hat{c} = CB$, then $\hat{b}$ has at least 10 letters and is rigid.

(b)  If $\hat{c} = BCBA \cdots CABC$, then it has at least 8 letters and is rigid.

(E)  The fifth case is  $ACBCA x ACBCA$ with matrix $\overline{BC} \cdots \overline{CB}$ (the shortest example is when there is just $\overline{BCB}$ and $x = BAC\ BAB\ CAB$).

(1) So if this is an instance of $\hat{a} \hat{b} \hat{a}$, then $x = \hat{b}$ has at least 9 letters and is rigid.

(2)  If this is an instance of $\hat{a}\hat{c}\hat{b}\hat{c}\hat{a}$, then $\hat{c} = B$ or $BACB \cdots BCAB$.

(a)  If $\hat{c} = B$, then $\hat{b}$ has at least 7 letters, which is disallowed by Theorem 1.

(b)  If $\hat{c} = BCBA \cdots CABC$, then it has at least 8 letters and is rigid.

(F)  The sixth case is  $ACBCA x ACBCA$ with matrix $\overline{CB} \cdots \overline{BC}$ (the shortest example is when there is just $\overline{CBC}$ and $x = CB\ ACA\ BC$).

(1) So if this is an instance of $\hat{a} \hat{b} \hat{a}$, then $x = \hat{b}$ has at least 7 letters, which is disallowed.

(2)  If this is an instance of $\hat{a}\hat{c}\hat{b}\hat{c}\hat{a}$, then $\hat{c} = C$ or $CBAC \cdots CABC$.

(a)  If $\hat{c} = C$, then $\hat{b} = BACAB$ or $BACA\cdots ACAB$.

(i)  If $\hat{b} = BACAB$, then $\hat{a}\hat{b}\hat{a}$ has a square on $CABA$.

(ii)  Otherwise, $\hat{b}$ has at least 8 letters and is rigid.

(b)  If $\hat{c} = CBAC \cdots CABC$, then it has at least 8 letters and is rigid.

\medskip

(II)  Next we study $\langle \hat{a}, \hat{b}, \hat{c} \rangle = \langle -,ACBCA,  - \rangle$ where two adjacent local occurrences of $\hat{b} = ACBCA$ are not rigid.  All but two cases give a minimum $x$ of length 14 or greater, so it divides with one keyword between $\hat{a}, \hat{c}$ of length at least 7, which is disallowed by Theorem 1.

(A)  The first case is  $ACBCA x ACBCA$ with matrix $\overline{BC} \cdots \overline{CB}$ (the shortest example is when there is just $\overline{BCB}$ and $x = BAC\ BAB\ CAB$).  No matter how the minimum or more general $x$ is divided, reversing the pieces gives an end-letter square.

(B)  The second case is  $ACBCA x ACBCA$ with matrix $\overline{CB} \cdots \overline{BC}$ (the shortest example is when there is just $\overline{CBC}$ and $x = CB\ ACA\ BC$).  No matter how $x$ is divided, reversing the pieces gives an end-letter square.

\medskip

(III)   Finally we address $\langle \hat{a}, \hat{b}, \hat{c} \rangle = \langle -,-,ACBCA \rangle$.  All six cases are covered by the length of the minimum $x$, which is at least 7.
\qed

\bigskip

There are six 4-letter subwords of $L$ up to rotation, and none of them is rigid.  Four of them have three matrices each, so the situation is complicated.  It is eased by having lengthy shortest gaps between many candidate keyword pairs.

\medskip

{\bf Theorem 7 (Four-letter keywords).}  {\it There is no local instance $[\hat{a},\hat{b},\hat{c},\rho]$ of $\kappa_2$ such that no keyword is locally rigid and the longest keyword has 4 letters.}

\medskip

{\it Proof.}  Sixty-four out of 112 major matrix-pair cases are resolved by having $x$s with lengths at least 14.  Twenty more are covered immediately by having $x$s with lengths at least 8.   We address the remaining 28 cases.

 The words we need to cover are listed, followed by their matrices and case loads.
\[
\begin{array}{lll}
ABAC &  (\overline{C}, \overline{C}) & \mbox{4 cases}\\
ABCA &(\overline{AB}, \overline{B}, \overline{C}) & \mbox{7 cases}\\ 
ABCB & (\overline{A}, \overline{A}) & \mbox{4 cases}\\
ACAB & (\overline{B}, \overline{AC},\overline{CA}) & \mbox{3 cases}\\
ACBA & (\overline{B}, \overline{BA},\overline{C}) & \mbox{7 cases}\\ 
ACBC & (\overline{A},\overline{BC}, \overline{CB}) & \mbox{3 cases} 
\end{array}
\]

{\bf ABAC}\ \ $(\overline{C}, \overline{C})$ \ \  4 cases 

(I)  We investigate $\langle \hat{a}, \hat{b}, \hat{c} \rangle = \langle ABAC, -, - \rangle$ where two adjacent local occurrences of $\hat{a} = ABAC$ are not rigid (because $\hat{a}$ is not locally rigid, such a pair exists).  The only matrix pair we need to consider is $\overline{C} \cdots \overline{C}$ (where the shortest example is just $\overline{C}$ and $x = BAB\ C$.  The pair that requires two $\overline{C}$s is covered by the excessive length of $x \geq 14$.

(A)  The first case is an instance of $\hat{a} \hat{b} \hat{a}$.  Therefore $x = \hat{b} = BABC$ for if $\hat{b}$ is longer than 4, this is contrary to assumption.  So $\hat{a}\hat{b}\hat{a}$ is rigid, and $\hat{c}$ is left flush.  Avoiding end-letter squares, $\hat{c}$ must be $ABC, BC$, or $BCAC$, lest it have five or more letters.  But all three have squares in $\hat{b}\hat{c}$.

(B)  The second case is an instance of $\hat{a}\hat{c}\hat{b}\hat{c}\hat{a}$.  Then $x = BABCA \cdots CBABC$ ($BABC$ is not an option).  So $\hat{c} = BABC$, for any longer is prohibited.  Then $\hat{b}$ with no more than 4 letters must be an initial segment of $ABAC$, so there is a square in $\hat{b}\hat{a}$.

(II)  The third case we study is $\langle \hat{a}, \hat{b}, \hat{c} \rangle = \langle -,ABAC,  - \rangle$ where two adjacent local occurrences of $\hat{b} = ABAC$ are not rigid.  The only matrix pair we need to consider is $\overline{C} \cdots \overline{C}$ (the shortest example is when there is just $\overline{C}$ and $x = BAB\ C$).  

(A)  If $x = BABC$, then splitting it into $\hat{a}$ and $\hat{c}$ gives one of $B, BA, BAB$, all of which are square with $\hat{b}$.  

(B)  If $x = BABC \cdots BABC$, then splitting it means one of $\hat{a}$ and $\hat{c}$ has at least 5 letters, which is prohibited.

(III)  We handle the fourth case $\langle \hat{a}, \hat{b}, \hat{c} \rangle = \langle -,-,ABAC \rangle$ by an argument similar to I.A.

\bigskip

{\bf ABCA}\ \ $(\overline{AB}, \overline{B},\overline{C})$ \ \ 7 cases

(I)  We investigate $\langle \hat{a}, \hat{b}, \hat{c} \rangle = \langle ABCA, -, - \rangle$ where two adjacent local occurrences of $\hat{a} = ABCA$ are not rigid.  

(A)  $ABCA x ABCA$ can have the matrix pattern $\overline{AB} \cdots \overline{B}$ (the shortest example is when there is just $\overline{AB}$ and $x = CBAC$). 

(1)  The first case is an instance of $\hat{a} \hat{b} \hat{a}$.  Then $x = \hat{b} = CBAC$ for if $\hat{b}$ cannot be longer than 4.  To avoid end-letter squares, $\hat{c}$ must be $B$ or have the form $B\cdots B$.  

(a)  If $\hat{c} = B$, then there is a square on $CBA$ in $\hat{b}\hat{c}\hat{a}$.

(b)  In the other case, $\hat{c} = BA \cdots B$ to avoid a square on $ABC$ in $\hat{a}\hat{c}$.  But now we see a square on $CBA$ in $\hat{b}\hat{c}$.

(2)  The second case is an instance of $\hat{a}\hat{c}\hat{b}\hat{c}\hat{a}$.  Therefore $x = \hat{c}\hat{b}\hat{c} =  CBAC\cdots CBAC$ ($CBAC$ is not an option).  So $\hat{c} = CBAC$, at maximal length.  This forces $\hat{b}$ to be an initial segment of $ABCA$, so there is a square in $\hat{b}\hat{a}$.

(B)  The third case is  $ABCA x ABCA$ with matrix $\overline{B}\cdots \overline{C}$ (the shortest example is when there is just $\overline{BC}$ and $x = CB\ CABAC\ B$).  The $\hat{a} \hat{b} \hat{a}$ possibility is handled by length, so we are left with $x = \hat{c}\hat{b}\hat{c}$.  The only solution is $\hat{c} = CB$ and $\hat{b} = CABA$, but that means there is a square in $\hat{b}\hat{a}$.

(C)  The fourth case is  $ABCA x ABCA$ with matrix $\overline{C}\cdots \overline{B}$ (the shortest example is when there is just $\overline{CB}$ and $x = BAC\ BCACB\ AC$).  The $\hat{a} \hat{b} \hat{a}$ possibility is handled by length, so we turn to $x = \hat{c}\hat{b}\hat{c}$.  The only solution is $\hat{c} = BAC$ and $\hat{b} = BCAC$, but that means there is a square on $ABC$ in $\hat{a}\hat{b}$.

\medskip

(II)  We investigate $\langle \hat{a}, \hat{b}, \hat{c} \rangle = \langle -, ABCA, - \rangle$ where two adjacent local occurrences of $\hat{b} = ABCA$ are not rigid.

(A)  The fifth case is  $ABCA x ABCA$ with matrix $\overline{AB} \cdots \overline{B}$ (the shortest example is when there is just $\overline{AB}$ and $x = CBAC$).  Splitting $x$ to obtain $\hat{c}$ and $\hat{a}$ and reversing the pieces gives a square on $C$.

(B)  The sixth case is  $ABCA x ABCA$ with matrix $\overline{B}\cdots \overline{C}$ (the shortest example is when there is just $\overline{BC}$ and $x = CB\ CABAC\ B$).  The only solution is to split the minimum $x$ in the middle, giving a square on $CB$ when reversed.

(C)  The remaining matrix pair $\overline{C}\cdots \overline{B}$ is covered by excessive length, since the shortest example of $x$ has 10 letters.

\medskip

(III)  We investigate $\langle \hat{a}, \hat{b}, \hat{c} \rangle = \langle -, -, ABCA \rangle$ where two adjacent local occurrences of $\hat{c} = ABCA$ are not rigid. 

(A)  The seventh case is  $ABCA x ABCA$ with matrix $\overline{AB} \cdots \overline{B}$ (the shortest example is when there is just $\overline{AB}$ and $x = CBAC$).  This is similar to I.A.1.

(B, C)  The other two matrix pairs are cleared via lengths of minimum $x$s.

\bigskip

{\bf ABCB}\ \ $(\overline{A}, \overline{A})$ \ \ 4 cases

This is ompletely similar to $ABAC$.

\bigskip

{\bf ACAB}\ \ $(\overline{B}, \overline{AC},\overline{CA})$ \ \ 3 cases

(I)  We investigate $\langle \hat{a}, \hat{b}, \hat{c} \rangle = \langle ACAB, -, - \rangle$ where two adjacent local occurrences of $\hat{a} = ACAB$ are not rigid.  

(A)  The first case is  $ACAB x ACAB$ with matrix $\overline{AC} \cdots \overline{CA}$ (the shortest example is when there is just $\overline{ACA}$ and $x = AC\ BAB\ CAB$).  The $\hat{a}\hat{b}\hat{a}$ possibility is handled by length, so we turn to $x = \hat{c}\hat{b}\hat{c}$.  There is no short solution for $\hat{c}$, so this case is complete.

(B)  The second case is  $ACAB x ACAB$ with matrix $\overline{CA} \cdots \overline{AC}$ (the shortest example is when there is just $\overline{CAC}$ and $x = CBA\ CBC\ ABCB$). The $\hat{a} \hat{b} \hat{a}$ possibility is handled by length, so we are left with $x = \hat{c}\hat{b}\hat{c}$.  The only feasible short solution for $\hat{c}$ is $CB$, which makes $\hat{b}$ at least 6 letters long.

(II)  We investigate $\langle \hat{a}, \hat{b}, \hat{c} \rangle = \langle -, ACAB, - \rangle$ where two adjacent local occurrences of $\hat{b} = ACAB$ are not rigid.  The third case, the only one not immediately covered by the minimum length of $x$, is  $ACAB x ACAB$ with matrix $\overline{AC} \cdots \overline{CA}$ (the shortest example is when there is just $\overline{ACA}$ and $x = AC\ BAB\ CAB$).  To keep keywords short, this minimum example must be divided into $ACBA$ and $BCAB$.  Both have end-letter squares with $\hat{b}$.

\bigskip

{\bf ACBA}\ \ $(\overline{B}, \overline{BA},\overline{C})$  \ \ 7 cases

(I)  We investigate $\langle \hat{a}, \hat{b}, \hat{c} \rangle = \langle ACBA, -, - \rangle$ where two adjacent local occurrences of $\hat{a} = ACBA$ are not rigid.

(A)  We start with $ACBA x ACBA$ having matrix $\overline{B}\cdots \overline{BA}$ (the shortest example is when there is just $\overline{BA}$ and $x = CABC$).  This has a development similar to $ABCA$ I.A.1 for the first case $\hat{b}$ and I.A.2 for second case $\hat{c}\hat{b}\hat{c}$.

(B)  The third case is $ACBA x ACBA$ with matrix $\overline{C}\cdots \overline{B}$ (the shortest example is when there is just $\overline{CB}$ and $x = B\ CABAC\ BC$) and can be immediately restricted to the $\hat{a}\hat{c}\hat{b}\hat{c}\hat{a}$ context.   The only solution avoiding immediate length disqualification is based on the minimum $x$ and has $\hat{c} = BC, \hat{b} = ABAC$.  The latter makes a square in $\hat{b}\hat{a}$.

(C) The fourth case is $ACBA x ACBA$ with matrix $\overline{B}\cdots \overline{C}$ (the shortest example is when there is just $\overline{BC}$ and $x = CA\ BCACB\ CAB$).  The only solution avoiding length disqualification is based on the minimum $x$ and has $\hat{c} = CAB, \hat{b} = CACB$.  The latter makes a square in $\hat{b}\hat{a}$.  Any other $\hat{c}$ will be longer than 4, no matter what $\hat{b}$ is.

(D)  For $ACBA x ACBA$ with the matrix pair $\overline{BA}\cdots \overline{B}$ is covered by having the minimum $x$ with 14 letters.

\medskip

(II)  Next we investigate $\langle \hat{a}, \hat{b}, \hat{c} \rangle = \langle -, ACBA, - \rangle$ where two adjacent local occurrences of $\hat{b} = ACBA$ are not rigid. 

(A)  The fifth case is  $ACBA x ACBA$ with matrix $\overline{B}\cdots \overline{BA}$ (the shortest example is when there is just $\overline{BA}$ and $x = CABC$).  So $x = CABC$ or $CABC \cdots CABC$, and however it is split into $\hat{a}$ and $\hat{c}$, reversing the pieces gives a square on $C$.

(B)  The sixth case is $ACBA x ACBA$ with matrix $\overline{C}\cdots \overline{B}$ (the shortest example is when there is just $\overline{CB}$ and $x = B\ CABAC\ BC$).  The only short split is on the minimum $x$:  $BCAB, ACBC$.  But these have a square on $BC$ when reversed.

(C, D) are covered by immediate length arguments.

\medskip

(III)  We investigate $\langle \hat{a}, \hat{b}, \hat{c} \rangle = \langle -, -, ACBA \rangle$ where two adjacent local occurrences of $\hat{c} = ACBA$ are not rigid. 

(A)  The seventh and only remaining case beyond quick length observations is $ACBA x ACBA$ with matrix $\overline{B}\cdots \overline{BA}$ (the shortest example is when there is just $\overline{BA}$ and $x = CABC$).  This is similar to I.A.$\hat{b}$, and also to $ABCA$ III.A.

(B, C, D) are covered by length arguments.

\bigskip

{\bf ACBC}\ \ $(\overline{A}, \overline{BC},\overline{CB})$ \ \ 3 cases

(I)  We investigate $\langle \hat{a}, \hat{b}, \hat{c} \rangle = \langle ACBC, -, - \rangle$ where two adjacent local occurrences of $\hat{a} = ACBC$ are not rigid (if $\hat{a}$ is not locally rigid, then such a pair exists).  

(A)  The first case not immediately covered by length arguments comes from $ACBC x ACBC$ with matrix $\overline{CB}\cdots \overline{BC}$ (the shortest example is when there is just $\overline{CBC}$ and $x = ACB\ ACA\ BC$).  We are only concerned with $x = \hat{c}\hat{b}\hat{c}$.  It is clear that $|\hat{c}| \geq 5$.

(B)  The second case is  $ACBC x ACBC$ with matrix $\overline{BC}\cdots \overline{CB}$ (the shortest example is when there is just $\overline{BCB}$ and $x = ABAC\ BAB\ CAB$). We are only concerned with $x = \hat{c}\hat{b}\hat{c}$.  It is clear that $|\hat{c}| \geq 5$.

(II)  Next we investigate $\langle \hat{a}, \hat{b}, \hat{c} \rangle = \langle -, ACBC, - \rangle$ where two adjacent local occurrences of $\hat{b} = ACBC$ are not rigid (if $\hat{b}$ is not locally rigid, then such a pair exists). 

(A)  The third case is $ACBC x ACBC$ with matrix $\overline{CB}\cdots \overline{BC}$ (the shortest example is when there is just $\overline{CBC}$ and $x = ACB\ ACA\ BC$).  The only short division into $\hat{a}$ and $\hat{c}$ is $ACBA, CABC$, using the minimum $x$.  Both words have end-letter squares with $\hat{a}$.

All remaining situations are immediately covered by length.    \qed

\bigskip

{\bf Theorem 8 (Short keywords).}  {\it There is no local instance $[\hat{a},\hat{b},\hat{c},\rho]$ of $\kappa_2$ such that no keyword is locally rigid and all the keywords have no more than 3 letters.}

\medskip

{\it Proof.}

(1) (Three letters maximum)

(a) One of the keywords has the form $XYX$.  The only other candidate forms are $ZYZ$ and $Z$, but not both at once, so three keywords are impossible.

(b)  One of the keywords has the form $XYZ$.  The only compatible forms are $XZ$ and $Y$.  This gives 12 cases up to rotation.

\medskip

$\langle ABC,AC,B\rangle:  ABC\ AC\ A\cdots $ has a square on $CA$.

$\langle ABC,B,AC\rangle:  ABC\ B\ ABC\ AC\ B \ AC \ ABC \ B \cdots \not \subterm L$.

$\langle AC,ABC,B\rangle:  AC\ ABC\ AC\ B\ ABC\ B\ AC\ A\cdots \not \subterm L$.

$\langle AC,B,ABC\rangle:  AC\ B\ AC\ ABC\ B\cdots \not \subterm L$.

$\langle B,AC,ABC\rangle:  B\ AC\ B\ ABC\ AC\cdots \not \subterm L$.

$\langle B,ABC,AC\rangle:  B\ ABC\ B\ AC\ A\cdots \not \subterm L$.

\medskip

$\langle ACB,AB,C\rangle:  ACB\ AB\ A\cdots $ has a square $BABA$.

$\langle ACB,C,AB\rangle:  ACB\ C\ ACB\ AB \cdots \not \subterm L$.

$\langle AB,ACB,C\rangle:  AB\ ACB\ AB\ C\ AC \cdots \not \subterm L$.

$\langle AB,C,ACB\rangle:  AB\ C\ AB\ ACB\ C\ ACB\ AB\cdots \not \subterm L$.

$\langle C,AB,ACB\rangle:  C\ AB\ C\ ACB\ AB\ A \cdots \not \subterm L$.

$\langle C,ACB,AB\rangle:  C\ ACB\ C\ AB\ ACB\ AB\ C\ AC \cdots \not \subterm L$.

\medskip

This last is the longest survivor, demonstrating that $\hat{a}\hat{b}\hat{a}\hat{c}\hat{b}\hat{c}\hat{a}\hat{b} \subterm \kappa_2$ is needed, as do $\langle ABC,B,AC\rangle$ and $\langle AC,ABC,B\rangle$.  The last $\hat{a}$ in the full $\kappa_2$ guarantees that we have both $\hat{a}\hat{b}\hat{a}$s when we use rigidity on both sides.

(2)  (Two letters maximum)  The longest keyword has the form $XY$.  The only compatible forms are $XZ, ZY, Z$, but no two at once, so three keywords are impossible.

(3)  (One letter)  The only solution is $X, Y, Z$, but $\kappa_2 = ABA\ CBC\ AB \cdots \not \subterm L$, and $\refmap(\kappa_2) = ACA\ BCB\ AC \cdots \not \subterm L$.
\qed

\bigskip

{\bf Theorem 9 (Word not in $L$).}  {\it There is no substitution instance of $\kappa_1$ (or $\kappa_3$ or $\kappa_4$) in $L$}.

\medskip

{\it Proof.}
The preceding results show that
\[
ABA\ CABCBAC\ B\ CABCBAC\ ABA
\]
and its rot variants are the only substitution instances of $\kappa_2$ in $L$, so the observation that
\[
ABA\ CABCBAC\ B\ CABCBAC\ ABA\ CA
\]
is already not a subterm of $L$ proves the absence of $\kappa_1$ from $L$.  That $\kappa_3$ is  also missing follows similarly.  Of course, $\kappa_4$ is also impossible.  \qed

\bigskip

Better ask bluntly about a shorter proof, a shorter word.

\end{document}